\documentclass[12pt]{article}
\usepackage{amsmath,amsfonts}
\usepackage{verbatim}
\usepackage{latexsym}
\usepackage{graphicx}
 \usepackage{color}
\makeatletter \@addtoreset{equation}{section}

\begin{document}

\newcommand{\E}{\mathbb{E}}
\newcommand{\PP}{\mathbb{P}}
\newcommand{\RR}{\mathbb{R}}

\newcommand{\Dt}{\D t}
\newcommand{\bX}{\bar X}
\newcommand{\bx}{\bar x}
\newcommand{\by}{\bar y}
\newcommand{\bp}{\bar p}
\newcommand{\bq}{\bar q}

\newtheorem{theorem}{Theorem}[section]
\newtheorem{lemma}[theorem]{Lemma}
\newtheorem{coro}[theorem]{Corollary}
\newtheorem{defn}[theorem]{Definition}
\newtheorem{assp}[theorem]{Assumption}
\newtheorem{expl}[theorem]{Example}
\newtheorem{prop}[theorem]{Proposition}
\newtheorem{rmk}[theorem]{Remark}

\newcommand\tq{{\scriptstyle{3\over 4 }\scriptstyle}}
\newcommand\qua{{\scriptstyle{1\over 4 }\scriptstyle}}
\newcommand\hf{{\textstyle{1\over 2 }\displaystyle}}
\newcommand\athird{{\scriptstyle{1\over 3 }\scriptstyle}}
\newcommand\hhf{{\scriptstyle{1\over 2 }\scriptstyle}}

\newcommand{\eproof}{\indent\vrule height6pt width4pt depth1pt\hfil\par\medbreak}

\def\a{\alpha} \def\g{\gamma}
\def\e{\varepsilon} \def\z{\zeta} \def\y{\eta} \def\o{\theta}
\def\vo{\vartheta} \def\k{\kappa} \def\l{\lambda} \def\m{\mu} \def\n{\nu}
\def\x{\xi}  \def\r{\rho} \def\s{\sigma}
\def\p{\phi} \def\f{\varphi}   \def\w{\omega}
\def\q{\surd} \def\i{\bot} \def\h{\forall} \def\j{\emptyset}

\def\be{\beta} \def\de{\delta} \def\up{\upsilon} \def\eq{\equiv}
\def\ve{\vee} \def\we{\wedge}

\def\F{{\cal F}}
\def\T{\tau} \def\G{\Gamma}  \def\D{\Delta} \def\O{\Theta} \def\L{\Lambda}
\def\X{\Xi} \def\S{\Sigma} \def\W{\Omega}
\def\M{\partial} \def\N{\nabla} \def\Ex{\exists} \def\K{\times}
\def\V{\bigvee} \def\U{\bigwedge}

\def\1{\oslash} \def\2{\oplus} \def\3{\otimes} \def\4{\ominus}
\def\5{\circ} \def\6{\odot} \def\7{\backslash} \def\8{\infty}
\def\9{\bigcap} \def\0{\bigcup} \def\+{\pm} \def\-{\mp}
\def\[{\langle} \def\]{\rangle}

\def\tl{\tilde}
\def\trace{\hbox{\rm trace}}
\def\diag{\hbox{\rm diag}}
\def\for{\quad\hbox{for }}
\def\refer{\hangindent=0.3in\hangafter=1}

\newcommand\wD{\widehat{\D}}

\thispagestyle{empty}

\title{
\bf The Truncated Euler--Maruyama Method for Stochastic Differential
Delay Equations}

\author{
{\bf Qian Guo${}^1$,
Xuerong Mao${}^2$\thanks{Corresponding author. E-mail: x.mao@strath.ac.uk.},
Rongxian Yue${}^1$}
\\
${}^1$ Department of Mathematics, \\
Shanghai Normal University, Shanghai, China. \\
${}^2$ Department of Mathematics and Statistics, \\
University of Strathclyde, Glasgow G1 1XH, U.K.
}
\date{}

\maketitle

\begin{abstract}

The numerical solutions of stochastic differential delay equations
(SDDEs)
under the generalized Khasminskii-type condition were discussed by Mao \cite{M11},
and the theory there showed that the Euler--Maruyama (EM) numerical solutions converge
to the true solutions \emph{in probability}. However,
there is so far no result on the strong convergence (namely in $L^p$)
of the numerical solutions for the SDDEs
under this generalized condition.
In this paper, we will use the truncated EM method developed by Mao \cite{M15} to study
 the strong convergence
of the numerical solutions for the SDDEs
under the generalized Khasminskii-type condition.

\medskip \noindent
{\small\bf Key words: } Brownian motion, stochastic differential
delay equation, It\^o's formula, truncated Euler--Maruyama,
Khasminskii-type condition.

\medskip \noindent
{\small\bf  Mathematical Subject Classifications  (2000)}:  60H10, 60J65.

\end{abstract}

\section{Introduction}

In the study of
 stochastic differential delay equations (SDDEs),
the classical existence-and-uniqueness theorem requires
the coefficients of the SDDEs satisfy the local Lipschitz condition
and the linear growth condition (see, e.g., \cite{DZ92,KM, M91, M94, Moh}).
However, there are many SDDEs which do not satisfy the
linear growth condition.  In 2002, Mao \cite{M02} generalized the
the well-known Khasminskii test \cite{Kha} from stochastic differential equations (SDEs) to SDDEs.   The Khasminskii-type theorem
established in \cite{M02} for SDDEs gives the
conditions, in terms of Lyapunov functions, under which the
solutions to SDDEs will not
explode to infinity at a finite time.  The Khasminskii-type theorem enables us to verify
if a given nonlinear SDDE has a unique global solution under the
local Lipschitz condition but without the linear growth condition.
In 2005, Mao and Rassias \cite{MR05}   demonstrated that
 there are many important SDDEs which are not covered by
the Khasminskii-type theorem given in  \cite{M02}, and established a
generalized Khasminskii-type theorem
which covers a very wide class of nonlinear SDDEs.

On the other hand, there are in general no explicit solutions to
nonlinear SDDEs, whence numerical solutions are required in
practice.  The numerical solutions under the linear growth condition
plus the local Lipschitz condition have  been
discussed intensively by many authors (see, e.g., \cite{BB,KloP,KP,M97, MS, MY06}).
The numerical solutions of SDDEs
under the generalized Khasminskii-type condition were discussed by Mao \cite{M11},
and the theory there showed that the Euler--Maruyama (EM) numerical solutions converge
to the true solutions \emph{in probability}. However,
there is so far no result on the strong convergence (namely in $L^p$)
of the numerical solutions for the SDDEs
under the generalized Khasminskii-type condition.

Recently, Mao \cite{M15} develops a new explicit numerical method,
called the truncated EM method, for SDEs under the Khasminskii-type condition
plus the local Lipschitz condition and establishes the strong convergence theory.
In this paper, we will use this new truncated EM method to study
 the strong convergence
of the numerical solutions for the SDDEs
under the generalized Khasminskii-type condition.

This paper is organized as follows: We will introduce necessary
notion, state the  generalized Khasminskii-type condition and
define the truncated EM numerical solutions for SDDEs in Section 2.  We will
establish the strong convergence theory for the truncated EM numerical solutions in
Sections 3 and 4 and discuss the convergence rates in Section 5.
In each of these three sections we will illustrate our theory
by examples. We will see from these examples that the truncated EM numerical method can
be applied to approximate the solutions of many highly nonlinear
SDDEs.  We will finally conclude our paper in Section 6.

\section{The Truncated Euler-Maruyama Method }

Throughout this paper, unless otherwise specified, we use the
following notation. Let $|\cdot|$ be the Euclidean norm in $\RR^n$. If
$A$ is a vector or matrix, its transpose is denoted by $A^T$. If $A$
is a matrix, its trace norm is denoted by $|A|=\sqrt{ \trace(A^TA)
}$. Let $\RR_+=[0,\8)$ and $\T >0$.  Denote by $C([-\T,0]; \RR^n)$   the
family of continuous functions from $[-\T,0]$ to $\RR^n$ with the norm
$\|\f\|=\sup_{-\T\le\o\le 0}|\f(\o)|$.
%Denote by $C(\RR^n\K [-\T,\8];
%\RR_+)$ the family of continuous functions from $\RR^n\K [-\T,\8]$ to
%$\RR_+$.
Let $(\W , {\cal F}, \{{\cal F}_t\}_{t\ge 0}, \PP)$ be a complete
probability space with a filtration  $\{{\cal F}_t\}_{t\ge 0}$
satisfying the usual conditions (i.e., it is increasing and right
continuous while ${\cal F}_0$ contains all $\PP$-null sets). Let
$B(t)= (B_1(t), \cdots, B_m(t))^T$ be an $m$-dimensional Brownian
motion defined on the probability space.
Moreover, for two real numbers $a$ and $b$, we use $a\ve b=\max(a,b)$
and $a\we b=\min(a,b)$. If $G$ is a set, its indicator function is denoted by
$I_G$, namely $I_G(x)=1$ if $x\in G$ and $0$ otherwise.
If $a$ is a real number, we denote by $\lfloor a\rfloor$ the largest integer
which is less or equal to $a$, e.g., $\lfloor -1.2\rfloor = -2$ and $\lfloor 2.3\rfloor =2$.

% Let $p>0$ and denote by  $L^p_{\F_t}([-\T,0]; \RR^n)$ the family of
% $\F_t$-measurable
% $C([-\T,0];\RR^n)$-valued random variables $\x=\{\x(\o):-\T\le\o\le 0\}$
%such that $\E\|\x\|^p<\8$.
%If $x(t)$ is an $\RR^n$-valued stochastic process on $t\in [-\T,\8)$,
%we let $x_t=\{x(t+\o):-\T\le \o \le 0\}$ for $t\ge 0$.

Consider a nonlinear SDDE
\begin{equation}
dx(t) = f(x(t),x(t-\T))dt+g(x(t), x(t-\T))dB(t), \quad t\ge 0,
\label{sdde}
\end{equation}
with the initial data  given by
\begin{equation}
\{x(\o):-\T\le\o\le 0\} =\x \in  C([-\T,0]; \RR^n). \label{id}
\end{equation}
Here
$$
f: \RR^n\K \RR^n  \to \RR^n  \quad\hbox{and}\quad
  g: \RR^n\K \RR^n  \to \RR^{n\K m}.
$$
We assume that the coefficients $f$ and $g$ obey the   Local
Lipschitz condition:

\begin{assp} \label{A2.1}
For every positive number $R$ there is a
positive constant $K_R$ such that
 $$
|f(x,y)-f(\bar x,\bar y)|^2 \ve |g(x,y)-g(\bar x,\bar y)|^2
 \le K_R(|x-\bar x|^2 + |y-\bar y|^2)
$$
for those $x,y,\bar x, \bar y\in \RR^n$ with $|x| \ve |y|\ve |\bar
x|\ve |\bar y| \le R$.
\end{assp}

The classical existence-and-uniqueness theorem does not only require
this local Lipschitz condition but also the linear growth condition
(see, e.g., \cite{M91, M94, M97, Moh}). In this paper we shall
retain the local Lipschitz condition but replace the linear growth
condition by a generalized Khasminskii-type condition.

\begin{assp} \label{A2.2}
There are constants $K_1>0$, $K_2\ge 0$ and $\be >2$ such that
\begin{equation}\label{KhasC}
x^T f(x,y) + \frac{1}{2} |g(x,y)|^2 \le K_1(1+|x|^2+|y|^2)
-K_2|x|^\be +K_2 |y|^\be
\end{equation}
for all $(x,y)\in \RR^n\K \RR^n$.
\end{assp}

To have a feeling about what type of nonlinear SDDEs to which our theory
may apply, please consider, for example,  the scalar SDDE
$$
dx(t) = [a_1 +a_2x^2(t-\T)- a_3x^3(t)]dt
+ [a_4|x(t)|^{3/2}+a_5|x(t-\T)|^{3/2}]dB(t), \quad t\ge 0,
$$
where $a_3>0$ and $a_1,a_2,a_4,a_5\in \RR$ (see Example \ref{E3.6}
for the details).
The following result,  established in \cite{MR05},  is a generalized
Khasminskii-type theorem on the existence and uniqueness of the
solution to the SDDE.

\begin{lemma} \label{L2.3}
Let Assumptions \ref{A2.1} and \ref{A2.2} hold. Then for any given
initial data (\ref{id}), there is a unique global solution $x(t)$ to
equation (\ref{sdde}) on $t \in [-\T,\8)$. Moreover, the solution
has the property that
\begin{equation}\label{bd}
\sup_{-\T\le t\le T} \E|x(t)|^2   < \8,  \quad\forall T>0.
\end{equation}
\end{lemma}

It has been shown (see, e.g., \cite{M11})
that under Assumptions \ref{A2.1} and \ref{A2.2},
the EM numerical solutions converge to the true solution in probability.
But, to our best knowledge, \emph{there is so far no result on the strong convergence under
these assumptions}.
In this paper, we will use the truncated EM method  developed in \cite{M15}
and show that the truncated EM solutions will converge to the true solution in $L^q$
for some $q\ge 1$.

To define the truncated EM numerical solutions, we first choose a
strictly increasing continuous
function $\mu:\RR_+\to\RR_+$ such that $\mu(r)\to\8$ as $r\to\8$ and
\begin{equation}\label{mudef}
\sup_{|x|\ve |y|\le r} \big( |f(x,y)|\ve |g(x,y)| \big) \le \mu(r), \quad \forall r\ge 1.
\end{equation}
Denote by $\mu^{-1}$ the inverse function of $\mu$ and we see that
$\mu^{-1}$ is a strictly increasing continuous
function from $[\mu(0),\8)$ to $\RR_+$.  We also choose
a constant $\D^*\in (0,1]$ and a strictly
decreasing function $h: (0,\D^*]\to (0,\8)$ such that
\begin{equation}\label{hdef}
h(\D^*)\ge \mu(1), \ \
\lim_{\D \to 0} h(\D) = \8 \quad\hbox{and}\quad
\D^{1/4} h(\D) \le 1, \ \ \forall \D \in (0,\D^*].
\end{equation}
For example, we may choose $\D^*\in (0,1)$ sufficiently small
such that $1/\D^* \ge (\m(1))^4$ and define $h(\D)=\D^{-1/4}$ for $\D\in (0,\D^*]$.
For a given step size $\D\in (0,\D^*]$, let us define a mapping
$\pi_\D$ from $\RR^n$ to the closed ball
$\{x\in\RR^n: |x|\le \mu^{-1}(h(\D))\}$ by
$$
\pi_\D(x) = (|x|\we \mu^{-1}(h(\D)))\, \frac{x}{|x|},
$$
where we set  $x/|x|=0$ when $x=0$.
That is, $\pi_\D$ will map $x$ to itself when $|x|\le \mu^{-1}(h(\D))$
and to $\mu^{-1}(h(\D)) x/|x|$ when $|x|> \mu^{-1}(h(\D))$.
We then define the truncated functions
\begin{equation}\label{truncated}
f_\D(x,y) = f(\pi_\D(x),\pi_\D(y))
\quad\hbox{and}\quad
g_\D(x,y) = g(\pi_\D(x),\pi_\D(y))
\end{equation}
for $x,y\in\RR^n$.
It is easy to see that
\begin{equation}\label{bdcoeff}
|f_\D(x,y)|\ve |g_\D(x,y)| \le \mu(\mu^{-1}(h(\D))) = h(\D),\quad
\forall x,y\in \RR^n.
\end{equation}
That is, both truncated functions $f_\D$ and $g_\D$ are bounded
although $f$ and $g$ may not.  More usefully, these truncated functions
preserve the generalized Khasminskii-type condition to a very nice degree as described in
the following lemma.

\begin{lemma}\label{L2.4}
 Let Assumption \ref{A2.2} hold.
Then, for every $\D\in (0,\D^*]$, we have
\begin{equation}\label{KhasC2}
x^T f_\D(x,y) + \frac{1}{2} |g_\D(x,y)|^2 \le 2K_1(1+|x|^2+|y|^2)
-K_2|\pi_\D(x)|^\be + K_2|\pi_\D(y)|^\be
\end{equation}
for all $x,y\in \RR^n$.
\end{lemma}

\noindent
{\it Proof}.  Fix any $\D \in (0,\D^*]$.  Recalling that
$h(\D^*)\ge \mu(1)$, we see that $\mu^{-1}(h(\D^*))\ge 1$. But
$h$ is decreasing while $\mu^{-1}$ is increasing, so
$\mu^{-1}(h(\D))\ge 1$.

For $x\in \RR^n$ with $|x|\le \mu^{-1}(h(\D))$ and any $y\in\RR^n$,
we have, by (\ref{KhasC}),
\begin{eqnarray}\label{2.10}
& & x^T f_\D(x,y) + \frac{1}{2} |g_\D(x,y)|^2 \nonumber \\
&=& \pi_\D(x)^T f(\pi_\D(x),\pi_\D(y)) + \frac{1}{2} |g(\pi_\D(x),\pi_\D(y))|^2
\nonumber \\
&\le & K_1(1+|\pi_\D(x)|^2+|\pi_\D(y)|^2) - K_2|\pi_\D(x)|^\be + K_2|\pi_\D(y)|^\be
\nonumber \\
&\le &  K_1(1+|x|^2+|y|^2) - K_2|\pi_\D(x)|^\be+ K_2|\pi_\D(y)|^\be,
\end{eqnarray}
which implies the desired assertion (\ref{KhasC2}).
On the other hand, for $x\in \RR^n$ with $|x|> \mu^{-1}(h(\D))$
and any $y\in\RR^n$, we have
\begin{eqnarray}\label{2.11}
& & x^T f_\D(x,y) + \frac{1}{2} |g_\D(x,y)|^2 \nonumber \\
&=& \pi_\D(x)^T f(\pi_\D(x),\pi_\D(y)) + \frac{1}{2} |g(\pi_\D(x),\pi_\D(y))|^2
\nonumber \\
&+& (x-\pi_\D(x))^T f(\pi_\D(x),\pi_\D(y))
\nonumber \\
&\le & K_1(1+|\pi_\D(x)|^2+|\pi_\D(y)|^2) - K_2|\pi_\D(x)|^\be
+K_2|\pi_\D(y)|^\be
\nonumber \\
&+ &  \Big( \frac{|x|}{\mu^{-1}(h(\D))}-1 \Big) \pi_\D(x)^T f(\pi_\D(x),\pi_\D(y)),
\end{eqnarray}
where (\ref{KhasC}) has been used. But once again we see from (\ref{KhasC})
that
\begin{eqnarray*}
& & \pi_\D(x)^T f(\pi_\D(x),\pi_\D(y)) \\
& \le & K_1(1+|\pi_\D(x)|^2+|\pi_\D(y)|^2) - K_2 [\mu^{-1}(h(\D))]^\be
+K_2 |\pi_\D(y)|^\be \\
&\le & K_1(1+|\pi_\D(x)|^2+|\pi_\D(y)|^2).
\end{eqnarray*}
Substituting this into (\ref{2.11}) yields
\begin{eqnarray}\label{2.12}
& & x^T f_\D(x,y) + \frac{1}{2} |g_\D(x,y)|^2 \nonumber \\
&\le & \frac{K_1|x|}{\mu^{-1}(h(\D))}(1+|\pi_\D(x)|^2+|\pi_\D(y)|^2)
- K_2 |\pi_\D(x)|^\be +K_2 |\pi_\D(y)|^\be \nonumber \\
&\le &  K_1|x|  (1+|x|+|y|)
- K_2 |\pi_\D(x)|^\be +K_2 |\pi_\D(y)|^\be \nonumber \\
&\le & 2K_1 (1+|x|^2+|y|^2)
- K_2 |\pi_\D(x)|^\be +K_2 |\pi_\D(y)|^\be.
\end{eqnarray}
Namely, we have showed that the required
assertion (\ref{KhasC2}) also holds for $x\in \RR^n$ with $|x|> \mu^{-1}(h(\D))$
and any $y\in\RR^n$. The proof is hence complete. $\Box$

From now on,  we will let the step size $\D$  be a \emph{fraction} of $\T$.
That is, we will use $\D=\T/M$ for some positive integer $M$.  When we use
the terms of a sufficiently
small $\D$, we mean that we choose $M$ sufficiently large.

Let us now form the discrete-time truncated EM solutions.
Define $t_k = k\D$ for $k=-M, -(M-1), \cdots, 0, 1, 2, \cdots$.
Set $X_\D(t_k) = \xi(t_k)$ for $k=-M, -(M-1), \cdots, 0$ and then form
\begin{equation} \label{TEM1}
 X_\D(t_{k+1})= X_\D(t_k)+f_\D(X_\D(t_k), X_\D(t_{k-M}))\D + g_\D(X_\D(t_k), X_\D(t_{k-M}))\D B_k
\end{equation}
for $k = 0, 1,2,\cdots$,  where $\D B_k=B(t_{k+1})-B(t_{k})$.  In our analysis, it is more convenient
to work on the continuous-time approximations.  There are two continuous-time
versions.  One is the continuous-time step process
$\bx_\D(t)$ on $t\in [-\T,\8)$ defined by
\begin{equation} \label{TEM2}
\bx_\D(t) = \sum_{k=-M}^\8 X_\D({t_k}) I_{[k\D, (k+1)\D)}(t).
\end{equation}
The other one is the continuous-time continuous process
$x_\D(t)$ on $t\in [-\T,\8)$ defined by $x_\D(t) = \xi(t)$ for $t\in [-\T,0]$
while for $t \ge 0$
\begin{equation} \label{TEM3a}
x_\D(t) =  \xi(0) + \int_0^t f_\D(\bx_\D(s),\bx_\D(s-\T))ds + \int_0^t g_\D(\bx_\D(s),\bx_\D(s-\T))dB(s).
\end{equation}
We see that $x_\D(t)$ is an It\^o process on $t\ge 0$ with its It\^o differential
\begin{equation} \label{TEM3}
dx_\D(t) =  f_\D(\bx_\D(t),\bx_\D(t-\T))dt +  g_\D(\bx_\D(t),\bx_\D(t-\T))dB(t).
\end{equation}
It is useful to know that $X_\D(t_k)=\bx_\D(t_k)=x_\D(t_k)$ for every $k\ge -M$,
namely they coincide at $t_k$.
Of course, $\bx_\D(t)$ is computable but $x_\D(t)$ is not in general.
However, the following lemma shows that $x_\D(t)$ and $\bx_\D(t)$ are close to each other
in the sense of $L^p$.  This indicates that it is sufficient to use $\bx_\D(t)$
in practice.  On the other hand, in our analysis, it is more convenient
to work on both of them.

\begin{lemma}\label{L2.5}
For any $\D \in (0,\D^*]$ and any $p\ge 2$, we have
\begin{equation}\label{2.17}
\E|x_\D(t)-\bar x_\D(t)|^p \le c_p \D^{p/2} (h(\D))^p, \quad \forall t\ge 0,
\end{equation}
where $c_p$ is a positive constant dependent only on $p$. Consequently
\begin{equation}\label{2.18}
\lim_{\D\to 0} \E|x_\D(t)-\bar x_\D(t)|^p =0, \quad \forall t\ge 0.
\end{equation}
\end{lemma}

\noindent
{\it Proof}. In what follows, we will use $c_p$ to stand for
generic positive real constants dependent only on $p$ and its values may change between occurrences.
Fix $\D \in (0,\D^*]$ arbitrarily.  For any $t\ge 0$, there is a unique integer
$k\ge 0$ such that $t_k\le t < t_{k+1}$. By (\ref{bdcoeff}) and
the properties of the It\^o integral (see, e.g., \cite{M97}), we then derive from (\ref{TEM3}) that
\begin{eqnarray*}
& & \E|x_\D(t)-\bar x_\D(t)|^p = \E|x_\D(t)- x_\D(t_k)|^p \\
&\le & c_p \Big( \E\Big|\int_{t_k}^t f_\D(\bx_\D(s),\bx_\D(s-\T)) ds \Big|^p
+ \E\Big|\int_{t_k}^t g_\D(\bx_\D(s),\bx_\D(s-\T)) dB(s) \Big|^p
\Big) \\
&\le & c_p \Big( \D^{p-1} \E \int_{t_k}^t |f_\D(\bx_\D(s),\bx_\D(s-\T))|^p ds
+ \D^{(p-2)/2} \E \int_{t_k}^t |g_\D(\bx_\D(s),\bx_\D(s-\T))|^p ds \Big) \\
&\le & c_p \D^{p/2} (h(\D))^p,
\end{eqnarray*}
which is (\ref{2.17}).
Noting from (\ref{hdef}) that
$\D^{p/2} (h(\D))^p\le \D^{p/4}$, we obtain (\ref{2.18}) from (\ref{2.17}) immediately.
$\Box$

\section{Convergence in $L^q$ for $q\in [1,2)$}

\emph{From now on we will fix $T>0$ arbitrarily}.
In this section we will show that
$$
\lim_{\D\to 0} \E|x_\D(T) - x(T)|^q = 0
\quad\hbox{and}\quad
\lim_{\D\to 0} \E|\bar x_\D(T) - x(T)|^q = 0
$$
for every $1\le q<2$.
By (\ref{bdcoeff}), it is obvious that for every $p\ge 2$,
$$
 \E|x_\D(t)|^p < \8, \quad\forall t\ge 0.
$$
The following lemma gives an upper bound, independent of $\D$, for the second moment.

\begin{lemma}\label{L3.1}
Let Assumptions \ref{A2.1} and \ref{A2.2} hold. Then
\begin{equation}\label{3.1}
\sup_{0<\D \le \D^*} \sup_{0\le t\le T} \E|x_\D(t)|^2 \le C,
\end{equation}
where, and from now on, $C$ stands for generic positive real constants dependent on
 $T, K_1, K_2,\xi$ (and $\bp$, $K_3$ etc. as well in the next sections)
but independent of $\D$ and its values may change between occurrences.
\end{lemma}

\noindent
{\it Proof}. Fix $\D \in (0,\D^*]$ and the initial data $\xi$ arbitrarily.
By the It\^o formula, we derive from (\ref{TEM3}) that for $0\le t \le T$,
\begin{eqnarray*}
\E|x_\D(t)|^2
& = & |\xi(0)|^2 +
\E\int_0^t  \Big( 2 x_\D^T(s) f_\D(\bx_\D(s),\bx_\D(s-\T))  +
 |g_\D(\bx_\D(s),\bx_\D(s-\T))|^2 \Big) ds \\
&= & |\xi(0)|^2 +
\E\int_0^t  \Big( 2 \bar x_\D^T(s) f_\D(\bx_\D(s),\bx_\D(s-\T)) +
  |g_\D(\bx_\D(s),\bx_\D(s-\T))|^2 \Big) ds \\
& + & \E\int_0^t  2(x_\D(s)- \bar x_\D(s))^T f_\D(\bx_\D(s),\bx_\D(s-\T)) ds.
\end{eqnarray*}
By Lemma \ref{L2.4},  we get
\begin{eqnarray} \label{3.2}
\E|x_\D(t)|^2
& \le & |\xi(0)|^2 +
4K_1 \E\int_0^t  (1+|\bx_\D(s)|^2+|\bx_\D(s-\T)|^2) ds
 \nonumber \\
& -& 2K_2 \E\int_0^t   |\pi_\D(\bx_\D(s))|^\be ds
+ 2K_2 \E\int_0^t |\pi_\D(\bx_\D(s-\T))|^\be  ds  \nonumber \\
& + & 2 \E\int_0^t |x_\D(s)- \bar x_\D(s)| | f_\D(\bx_\D(s),\bx_\D(s-\T))| ds.
\end{eqnarray}
However, it is easy to show that
\begin{eqnarray} \label{3.3}
& &  |\xi(0)|^2 +
4K_1 \E\int_0^t  (1+|\bx_\D(s)|^2+|\bx_\D(s-\T)|^2) ds \nonumber \\
& & \quad
\le C + 8K_1 \int_0^t \Big( \sup_{0\le u\le s} \E|x_\D(u)|^2 \Big)ds.
\end{eqnarray}
Moreover,
\begin{eqnarray} \label{3.4}
& & -2K_2 \E\int_0^t  |\pi_\D(\bx_\D(s))|^\be ds
+2K_2 \E\int_0^t  |\pi_\D(\bx_\D(s-\T))|^\be  ds  \nonumber \\
&= & -2K_2 \E\int_0^t |\pi_\D(\bx_\D(s))|^\be ds
+ 2K_2 \E\int_{-\T}^{t-\T} |\pi_\D(\bx_\D(s))|^\be ds \nonumber \\
& \le & 2K_2 \int_{-\T}^{0} |\pi_\D(\bx_\D(s))|^\be ds   \le  2\T K_2 \|\xi\|^\be.
\end{eqnarray}
Furthermore,  by Lemma \ref{L2.5} with $p=2$ and inequalities (\ref{bdcoeff}) and (\ref{hdef}),
we derive that
\begin{eqnarray}\label{3.5}
& & \E\int_0^t |x_\D(s)- \bar x_\D(s)| | f_\D(\bx_\D(s),\bx_\D(s-\T))|ds
\nonumber \\
&\le & h(\D) \int_0^T \E|x_\D(s)- \bar x_\D(s)| ds
\nonumber \\
&\le & h(\D)\int_0^T (\E|x_\D(s)- \bar x_\D(s)|^{2}) ^{1/2} ds
\nonumber \\
&\le & C  (h(\D))^{2}  \D^{1/2} \le C.
\end{eqnarray}
Substituting (\ref{3.3})-(\ref{3.5}) into (\ref{3.2}) yields
\begin{eqnarray*}
\E|x_\D(t)|^2   \le   C + 8K_1 \int_0^t \Big( \sup_{0\le u\le s} \E|x_\D(u)|^2 \Big)ds.
\end{eqnarray*}
As this holds for any $t\in [0,T]$ while the sum of the right-hand-side (RHS) terms
is non-decreasing in $t$, we then see
\begin{eqnarray*}
\sup_{0\le u\le t} \E|x_\D(u)|^2
\le C + 8K_1 \int_0^t \Big( \sup_{0\le u\le s} \E|x_\D(u)|^2 \Big)ds.
\end{eqnarray*}
The well-known Gronwall inequality yields that
$$
\sup_{0\le u\le T} \E|x_\D(u)|^2 \le C.
$$
As this holds for any $\D \in (0,\D^*]$ while $C$ is independent of $\D$, we obtain the
required assertion (\ref{3.1}). $\Box$

Let us present two more lemmas before we state one of our main results in this paper.

\begin{lemma}\label{L3.2}
Let Assumptions \ref{A2.1} and \ref{A2.2} hold. For any real number $R> \|\xi\|$,
define the stopping time
$$
\T_R=\inf\{t\ge 0: |x(t)|\ge R\},
$$
where throughout this paper we set $\inf\emptyset = \8$ (and as usual
$\emptyset$ denotes the empty set). Then
\begin{equation}\label{3.6}
\PP(\T_R\le T) \le \frac{C}{R^2}.
\end{equation}
(Recall that $C$ stands for generic positive real constants dependent on $T,K_1,K_2,\xi$
so $C$ here is independent of $R$.)
\end{lemma}

\noindent
{\it Proof}. By the It\^o formula and Assumption \ref{A2.2}, we derive that
for $0\le t\le T$,
\begin{eqnarray*}
\E|x(t\we\T_R)|^2 &\le& |\xi(0)|^2 + 2K_1\E\int_0^{t\we\T_R} (1+|x(s)|^2+|x(s-\T)|^2 ) ds
\\
&- &2K_2 \E\int_0^{t\we\T_R} |x(s)|^\be ds
+ 2K_2 \E\int_0^{t\we\T_R} |x(s-\T)|^\be ds
\\
&\le & |\xi(0)|^2 + 2K_1T + 2K_1\E  \int_0^t \big( |x(s\we\T_R)|^2
  + |x((s-\T)\we\T_R)|^2 \big) ds  \\
&+ & 2K_2 \int_{-\T}^{0} |\xi(s)|^\be  ds \\
&\le & C+  2K_1  \int_0^t \big( \E|x(s\we\T_R)|^2
  + \E|x((s-\T)\we\T_R)|^2 \big) ds \\
  &\le & C + 4K_1  \int_0^t \Big( \sup_{0\le u\le s}\E|x(u\we\T_R)|^2
   \Big) ds
\end{eqnarray*}
But the sum of the RHS terms is non-decreasing in $t$, we hence have
$$
\sup_{0\le u\le t}\E|x(u\we\T_R)|^2
\le C + 4K_1  \int_0^t \Big( \sup_{0\le u\le s}\E|x(u\we\T_R)|^2
   \Big) ds.
$$
The Gronwall inequality shows
$$
 \sup_{0\le u\le T}\E|x(u\we\T_R)|^2 \le C.
$$
In particular, we have
$$
\E|x(T\we\T_R)|^2 \le C.
$$
This implies, by the Chebyshev inequality,
$$
R^2\, \PP(\T_R\le T) \le C
$$
and the assertion follows. $\Box$

\begin{lemma}\label{L3.3}
Let Assumptions \ref{A2.1} and \ref{A2.2} hold. For any real number $R> \|\xi\|$
and $\D\in (0,\D^*]$, define the stopping time
$$
\r_{\D,R} =\inf\{t\ge 0: |x_\D(t)|\ge R\}.
$$
Then
\begin{equation}\label{3.7}
\PP(\r_{\D,R}\le T) \le \frac{C}{R^2}.
\end{equation}
(Please recall that $C$ is independent of $\D$ and $R$.)
\end{lemma}

\noindent
{\it Proof}. We simply write $\r_{\D,R}= \r$.  In the same way as (\ref{3.2})
was obtained, we can show that
for $0\le t\le T$,
\begin{eqnarray} \label{3.8}
\E|x_\D(t\we\r)|^2
& \le & |\xi(0)|^2 +
4K_1 \E\int_0^{t\we\r}  (1+|\bx_\D(s)|^2+|\bx_\D(s-\T)|^2) ds
 \nonumber \\
& -& 2K_2 \E\int_0^{t\we\r}   |\pi_\D(\bx_\D(s))|^\be ds
+ 2K_2 \E\int_0^{t\we\r} |\pi_\D(\bx_\D(s-\T))|^\be  ds  \nonumber \\
& + & 2 \E\int_0^{t\we\r} |x_\D(s)- \bar x_\D(s)| | f_\D(\bx_\D(s),\bx_\D(s-\T))| ds.
\end{eqnarray}
In the same way as we performed in the proofs of Lemmas \ref{L3.1} and \ref{L3.2}, we can then show that
\begin{eqnarray} \label{3.9}
\E|x_\D(t\we\r)|^2
& \le & C +
8K_1 \int_0^t  \Big( \sup_{0\le u \le s} \E|\bx_\D(u\we\r)|^2 \Big) ds
 \nonumber \\
& + & 2 \E\int_0^t |x_\D(s)- \bar x_\D(s)| | f_\D(\bx_\D(s),\bx_\D(s-\T))| ds.
\end{eqnarray}
This, together with (\ref{3.5}), implies
$$
\E|x_\D(t\we\r)|^2
\le C +
8K_1 \int_0^t  \Big( \sup_{0\le u \le s} \E |\bx_\D(u\we\r)|^2 \Big) ds.
$$
Noting that the sum of the RHS terms is increasing in $t$ while
$$
\sup_{0\le u \le s} \E |\bx_\D(u\we\r)|^2 \le \sup_{0\le u \le s} \E |x_\D(u\we\r)|^2,
$$
we get
$$
 \sup_{0\le u \le t} \E|x_\D(u\we\r)|^2
\le C +
8K_1 \int_0^t  \Big( \sup_{0\le u \le s} \E |x_\D(u\we\r)|^2 \Big) ds.
$$
The Gronwall inequality shows
$$
 \sup_{0\le u \le T} \E|x_\D(u\we\r)|^2 \le C.
$$
This implies the required assertion (\ref{3.7}) easily.  $\Box$

For the numerical solutions to converge to the true solution in $L^q$,
we need to assume that the initial data are H\"older
continuous with exponent $\g$ (or $\g$-H\"older continuous).  This is
a standard condition which is also needed for the classical EM method under the
global Lipschitz condition (see, e.g.,
\cite{MS,MY06,WM08}).

\begin{assp}\label{A5.2}
There is a pair of  constants $K_3>0$ and $\g\in (0,1]$ such that the initial
data $\xi$ satisfies
$$
|\xi(u)-\xi(v)|\le K_3|u-v|^\g, \quad -\T\le v<u\le 0.
$$
\end{assp}

We can now show one of our main results in this paper.

\begin{theorem}\label{T3.4}
Let Assumptions \ref{A2.1}, \ref{A2.2} and \ref{A5.2}  hold.  Then, for any $q\in [1,2)$,
\begin{equation}\label{3.10}
\lim_{\D\to 0} \E|x_\D(T) - x(T)|^q = 0
\quad\hbox{and}\quad
\lim_{\D\to 0} \E|\bar x_\D(T) - x(T)|^q = 0.
\end{equation}
\end{theorem}

\noindent
{\it Proof}.   Let $\T_R$ and $\r_{\D,R}$ be the same as before. Set
$$
\o_{\D,R}=\T_R \we\r_{\D,R}
\quad\hbox{and} \quad
e_\D(T)= x_\D(T) - x(T).
$$
Obviously
\begin{equation}\label{3.11}
\E|e_\D(T)|^q = \E\Big( |e_\D(T)|^q I_{\{\o_{\D,R}>T\}} \Big)
+ \E\Big( |e_\D(T)|^q I_{\{\o_{\D,R}\le T\} } \Big).
\end{equation}
Let $\de >0$ be arbitrary.  Using the Young inequality
$$
a^q b = (\de a^2)^{q/2} \Big( \frac{b^{2/(2-q)}}{\de^{q/(2-q)} }\Big)^{(2-q)/2}
\le \frac{q\de}{2} a^2 + \frac{2-q}{ 2\de^{q/(2-q)} } b^{2/(2-q)}, \quad\forall a,b>0,
$$
we have
$$
\E\Big( |e_\D(T)|^q I_{\{\o_{\D,R}\le T\} } \Big) \le
\frac{q\de}{2} \E|e_\D(T)|^2 + \frac{2-q}{2\de^{q/(2-q)}} \PP(\o_{\D,R}\le T).
$$
By Lemmas \ref{L2.3} and \ref{L3.1}, we have
$$
\E|e_\D(T)|^2 \le C,
$$
while by Lemmas \ref{L3.2} and \ref{L3.3},
$$
\PP(\o_{\D,R}\le T) \le \PP(\T_R\le T)+\PP(\r_{\D,R}\le T) \le \frac{C}{R^2}.
$$
We hence have
$$
\E\Big( |e_\D(T)|^q I_{\{\o_{\D,R}\le T\} } \Big)
\le \frac{C q\de }{2} + \frac{C(2-q)}{2 R^2 \de^{q/(2-q)}}.
$$
Substituting this into (\ref{3.11}) yields
\begin{equation} \label{3.12}
\E|e_\D(T)|^q \le \E\Big( |e_\D(T )|^q I_{\{\o_{\D,R}>T\}} \Big)
+ \frac{C q\de }{2} + \frac{C(2-q)}{2 R^2 \de^{q/(2-q)}}.
\end{equation}
Now, let $\e>0$ be arbitrary. Choose $\de$ sufficiently small for
$Cq\de/2 \le \e/3$ and then choose $R$ sufficiently large for
$$
\frac{C(2-q)}{2 R^2 \de^{q/(2-q)}} \le \frac{\e}{3}.
$$
We then see from (\ref{3.12}) that for this particularly chosen $R$,
\begin{equation} \label{3.13}
\E|e_\D(T)|^q \le \E\Big( |e_\D(T )|^q I_{\{\o_{\D,R}>T\}} \Big)
+\frac{2\e}{3}.
\end{equation}
If we can show that for all sufficiently small $\D$,
\begin{eqnarray} \label{3.14}
\E\Big( |e_\D(T )|^q I_{\{\o_{\D,R}>T\}} \Big) \le \frac{ \e}{3},
\end{eqnarray}
we have
$$
\lim_{\D\to 0} \E|e_\D(T)|^q = 0,
$$
and then by Lemma \ref{L2.5}, we also have
$$
\lim_{\D\to 0} \E|x(T)-\bar x_\D(T)|^q = 0.
$$
In other words, to complete our proof, all we need is to show (\ref{3.14}). For this purpose, we define the truncated functions
$$
F_R(x,y) = f\Big( (|x|\we R) \frac{x}{|x|},(|y|\we R) \frac{y}{|y|} \Big)
\quad\hbox{and}\quad
G_R(x,y) = g\Big( (|x|\we R) \frac{x}{|x|} ,(|y|\we R) \frac{y}{|y|}\Big)
$$
for $x,y\in\RR^n$.
Without loss of any generality, we may assume that $\D^*$ is already
sufficiently small for $\mu^{-1}(h(\D^*)) \ge R$. Hence, for all $\D\in (0,\D^*]$,
we have that
$$
f_\D(x,y) = F_R(x,y)
\quad\hbox{and}\quad
g_\D(x,y) = G_R(x,y)
$$
for those $x,y\in\RR^n$ with $|x|\ve|y|\le R$.
Consider the SDDE
\begin{equation} \label{sdde2}
dz(t) = F_R(z(t),z(t-\T))dt + G_R(z(t),z(t-\T))dB(t)
\end{equation}
on $t\ge 0$ with the initial data $z(u)=\xi(u)$ on $u\in [-\T,0]$.
By Assumption \ref{A2.1}, we see that both $F_R(x,y)$ and $G_R(x,y)$ are
globally Lipschitz continuous with the Lipschitz constant $K_R$. So the SDDE
(\ref{sdde2}) has a unique global solution $z(t)$ on $t\ge -\T$. It is straightforward
to see that
\begin{equation} \label{3.16}
\PP\{x(t\we\T_R) = z(t\we\T_R) \hbox{ for all } 0\le t\le T\} = 1.
\end{equation}
On the other hand, for each step size $\D\in (0,\D^*]$, we can apply the (classical) EM method
to the SDDE (\ref{sdde2}) and we denote by $z_\D(t)$ the continuous-time continuous EM solution.
It is again straightforward
to see that
\begin{equation} \label{3.17}
\PP\{x_\D(t\we\r_{\D,R}) = z_\D(t\we\r_{\D,R})  \hbox{ for all } 0\le t\le T\} = 1.
\end{equation}
However, it is well known (see, e.g., \cite{MS,MY06}) that
\begin{equation} \label{3.18}
\E\Big( \sup_{0\le t\le T} |z(t)-z_\D(t)|^q \Big) \le H \D^{q(0.5\we\g)},
\end{equation}
where $H$ is a positive constant dependent on $K_R, T, \xi, q$ but independent
of $\D$.
Consequently,
$$
\E\Big( \sup_{0\le t\le T} |z(t\we\o_{\D,R})-z_\D(t\we\o_{\D,R})|^q \Big) \le H \D^{q(0.5\we\g)}.
$$
Using (\ref{3.16}) and (\ref{3.17}), we then have
\begin{equation}\label{3.19}
\E\Big( \sup_{0\le t\le T} |x(t\we\o_{\D,R})-x_\D(t\we\o_{\D,R})|^q \Big) \le H \D^{q(0.5\we\g)},
\end{equation}
which implies
$$
\E\Big( |x(T\we\o_{\D,R})-x_\D(T\we\o_{\D,R})|^q \Big) \le H \D^{q(0.5\we\g)}.
$$
Finally
\begin{eqnarray} \label{3.20}
& & \E\Big( |e_\D(T )|^q I_{\{\o_{\D,R}>T\}} \Big)
= \E\Big( |e_\D(T\we\o_{\D,R} )|^q I_{\{\o_{\D,R}>T\}} \Big) \nonumber \\
& &\quad \le \E\Big( |x(T\we\o_{\D,R})-x_\D(T\we\o_{\D,R})|^q \Big)
\le H \D^{q(0.5\we\g)}.
\end{eqnarray}
This implies (\ref{3.14}) as desired.
The proof is therefore complete. $\Box$

Let make a useful remark which will be used in next sections before we
discuss an example to illustrate our theory.

\begin{rmk} \label{R3.5}
It is known (see, e.g., \cite{MS,MY06}) that (\ref{3.18}) holds for any $q\ge 2$.
We hence see from the proof above that both (\ref{3.19}) and (\ref{3.20})
hold for any $q\ge 2$ too.
\end{rmk}

\begin{expl} \label{E3.6}
{\rm  Consider the scalar SDDE
\begin{equation} \label{E6.1a}
dx(t) = x(t)\Big([a_1 +a_2x(t-\T)- a_3x^2(t)]dt + [a_4x(t)+a_5x(t-\T)]dB(t)\Big), \quad t\ge 0,
\end{equation}
with the initial data $\{x(\o):-\T\le\o\le 0\}=\xi\in C([-\T,0];(0,\8))$, where $B(t)$ is a
scalar Brownian motion and
$a_i$ ($1\le i\le 5$) are all positive numbers with
\begin{equation} \label{E6.1b}
a_3 > a_4^2+a_5^2.
\end{equation}
This is a stochastic delay population system (see, e.g., \cite{BM04a,BM04b,MYZ}).
It can be shown that given the initial data $\{x(\o):-\T\le\o\le 0\}=\xi\in C([-\T,0];(0,\8))$,
the solution will remain positive for all $t\ge 0$ with probability 1.   We can therefore
regard equation (\ref{E6.1a}) as an SDDE in $\RR$ with the coefficients
$$
f(x,y) = x(a_1 +a_2y- a_3x^2)  \quad\hbox{and}\quad
g(x,y) = x(a_4x+ a_5y), \quad x,y\in\RR.
$$
It is obvious that these coefficients are locally Lipschitz continuous, namely, they satisfy Assumption \ref{A2.1}.
We also assume that the initial data satisfy Assumption \ref{A5.2}.
Moreover,  we set $\de=  a_3 - a_4^2-a_5^2$, which is positive by (\ref{E6.1b}),  and derive
\begin{eqnarray*}
xf(x,y)+\frac{1}{2} |g(x,y)|^2
&\le& a_1x^2+a_2x^2|y|-a_3x^4+a_4^2x^4+a_5^2x^2y^2 \\
& \le & a_1x^2+(a_2^2/4\de) y^2-(a_3-\de-a_4^2-0.5a_5^2)x^4 + 0.5a_5^2y^4 \\
&\le &  (a_1\ve (a_2^2/4\de))(1+x^2+y^2) -0.5a_5^2x^4+ 0.5a_5^2 y^4.
\end{eqnarray*}
That is, Assumption \ref{A2.2} is satisfied as well.  We can
therefore apply the truncated EM method to obtain the
numerical solutions of the SDDE (\ref{E6.1a}).
For this purpose, we observe that, for $r\ge 1$,
$$
\sup_{|x|\ve|y|\le r} (|f(x,y)|\ve |g(x,y)|)
\le (a_1 r+a_2r^2+a_3r^3)\ve ((a_4+a_5)r^2) \le a r^3,
$$
where $a=(a_1+a_2+a_3)\ve(a_4+a_5)$.   We can therefore define
$\mu:\RR_+\to\RR_+$ by
$$
\mu(r) = a r^3, \quad r\ge 0.
$$
Its inverse function $\mu^{-1}: \RR_+ \to \RR_+$ has the form
$$
\mu^{-1}(r) = \Big( \frac{r}{a}  \Big)^{1/3}, \quad r\ge 0.
$$
Let $\r\in (0,1/4]$ and $\D^*=(1\ve(8a))^{-1/\r} \in (0,1]$.  Define
$h(\D)=\D^{-\r}$ for $\D\in (0,\D^*]$.  We then see that
$h(\D^*)\ge 8a=\mu(2)$, $\lim_{\D\to 0} h(\D)=\8$ and
$$
\D^{1/4} h(\D) = \D^{1/4-\r} \le 1, \quad \forall \D\in (0,\D^*]
$$
as required by (\ref{hdef}).   With these chosen functions $\mu$ and $h$,
we can then apply the truncated EM method to obtain the numerical
solutions $x_\D(t)$ and $\bx_\D(t)$ of the SDDE (\ref{E6.1a}).
Moreover, Theorem \ref{T3.4} shows that these numerical solutions
will converge to the true solution $x(t)$  in the sense that
$$
\lim_{\D\to 0} \E|x_\D(t) - x(t)|^q = 0
\quad\hbox{and}\quad
\lim_{\D\to 0} \E|\bar x_\D(t) - x(t)|^q = 0
$$
for any $q\in [1,2)$.
}
\end{expl}

\section{Convergence in $L^q$ for $q\ge 2$}

In the previous section, we showed that the truncated EM solutions
$x_\D(T)$ and $\bar x_\D(T)$ will converge to the true solution $x(T)$
in $L^q$ for any $q\in [1,2)$. This is sufficient for
some applications, for example,
when we need to approximate the mean value of the solution
or the European call option value (see, e.g., \cite{HM2005}).
However, we sometimes need to approximate
the variance or higher moment of the solution.
In these situations, we need to have the convergence
in $L^q$ for $q\ge 2$. For this purpose, we impose
a stronger Khasminskii-type condition.

\begin{assp} \label{A4.1}
There is a pair of constants $\bp>2$ and $K_1>0$ such that
\begin{equation}\label{KhasC3}
x^T f(x,y) + \frac{\bp-1}{2} |g(x,y)|^2 \le K_1(1+|x|^2+|y|^2)
\end{equation}
for all $(x,y)\in \RR^n\K \RR^n$.
\end{assp}

Once again, the truncated functions $f_\D$ and $g_\D$ preserve
this condition nicely.

\begin{lemma}\label{L4.2}
Let Assumption \ref{A4.1} hold.  Then, for every $\D\in (0,\D^*]$, we have
\begin{equation}\label{KhasC4}
x^T f_\D(x,y) + \frac{\bp-1}{2} |g_\D(x,y)|^2 \le 2K_1(1+|x|^2+|y|^2)
\end{equation}
for all $x,y\in \RR^n$.
\end{lemma}

This lemma can be proved in the same way as Lemma \ref{L2.4} was proved.
We also cite a stronger result than Lemma \ref{L2.3} from \cite{MR05}.

\begin{lemma} \label{L4.3}
Let Assumptions \ref{A2.1} and \ref{A4.1} hold. Then for any given
initial data (\ref{id}), there is a unique global solution $x(t)$ to
equation (\ref{sdde}) on $t \in [-\T,\8)$. Moreover, the solution
has the property that
\begin{equation}\label{bd2}
\sup_{-\T\le t\le T} \E|x(t)|^{\bp}  < \8.
\end{equation}
\end{lemma}

Let us now establish a stronger result than Lemma \ref{L3.1}.

\begin{lemma}\label{L4.4}
Let Assumptions \ref{A2.1} and \ref{A4.1} hold. Then
\begin{equation}\label{4.4}
\sup_{0<\D \le \D^*} \sup_{0\le t\le T} \E|x_\D(t)|^{\bp} \le C.
\end{equation}
\end{lemma}

\noindent
{\it Proof}. Fix any $\D \in (0,\D^*]$. By the It\^o formula, we derive from (\ref{TEM3}) that, for $0\le t \le T$,
\begin{eqnarray*}
\E|x_\D(t)|^{\bp}
& \le & |\xi(0)|^{\bp} +
\E\int_0^t \bp |x_\D(s)|^{\bp-2} \\
& & \quad \K \Big( x_\D^T(s) f_\D(\bx_\D(s),\bx_\D(s-\T))
 +
\frac{\bp-1}{2} |g_\D(\bx_\D(s),\bx_\D(s-\T))|^2 \Big) ds \\
&= & |\xi(0)|^{\bp} +
\E\int_0^t \bp |x_\D(s)|^{\bp-2}
\\
& & \quad \K \Big( \bar x_\D^T(s) f_\D(\bx_\D(s),\bx_\D(s-\T)) +
\frac{\bp-1}{2} |g_\D(\bx_\D(s),\bx_\D(s-\T))|^2 \Big) ds \\
& + & \E\int_0^t \bp |x_\D(s)|^{\bp-2} (x_\D(s)- \bar x_\D(s))^T f_\D(\bx_\D(s),\bx_\D(s-\T)) ds.
\end{eqnarray*}
By Lemma \ref{L4.2} and the Young inequality
$$
a^{\bp-2}b \le \frac{\bp-2}{\bp} \, a^{\bp} +\frac{2}{\bp}\, b^{\bp/2}, \quad\forall a, b\ge 0,
$$
we then have
\begin{eqnarray*}
\E|x_\D(t)|^{\bp}
& \le & |\xi(0)|^{\bp} +
\E\int_0^t  2 \bp K_1 |x_\D(s)|^{\bp-2} (1+|\bx_\D(s)|^2+|\bx_\D(s-\T)|^2) ds \\
& + & (\bp-2) \E\int_0^t |x_\D(s)|^{\bp} ds  \\
&+ & 2 \E\int_0^t |x_\D(s)- \bar x_\D(s)|^{\bp/2} | f_\D(\bx_\D(s),\bx_\D(s-\T))|^{\bp/2} ds \\
& \le & C+ C
\int_0^t \big( \E|x_\D(s)|^{\bp} + \E|\bar x_\D(s)|^{\bp}
+ \E|\bar x_\D(s-\T)|^{\bp} \big) ds \\
& + & 2 \E\int_0^T |x_\D(s)- \bar x_\D(s)|^{\bp/2}
 | f_\D(\bx_\D(s),\bx_\D(s-\T))|^{\bp/2} ds.
\end{eqnarray*}
But, by Lemma \ref{L2.5} with $p=\bp$ and inequalities (\ref{bdcoeff}) and (\ref{hdef}),
we have
\begin{eqnarray}\label{4.5}
& &\E\int_0^T |x_\D(s)- \bar x_\D(s)|^{\bp/2} | f_\D(\bx_\D(s),\bx_\D(s-\T))|^{\bp/2} ds
\nonumber \\
&\le & (h(\D))^{\bp/2} \int_0^T \E(|x_\D(s)- \bar x_\D(s)|^{\bp/2}) ds
\nonumber \\
&\le & (h(\D))^{\bp/2} \int_0^T (\E|x_\D(s)- \bar x_\D(s)|^{\bp}) ^{1/2} ds
\nonumber \\
&\le & c_{\bp} T (h(\D))^{\bp} \D^{\bp/4} \le c_{\bp} T.
\end{eqnarray}
We therefore have
\begin{eqnarray*}
\E|x_\D(t)|^{\bp}
& \le & C+ C \int_0^t \big( \E|x_\D(s)|^{\bp} + \E|\bar x_\D(s)|^{\bp} + \E|\bar x_\D(s-\T)|^{\bp} \big) ds \\
& \le & C + C \int_0^t \Big( \sup_{0\le u\le s} \E|x_\D(u)|^{\bp} \Big)ds.
\end{eqnarray*}
As this holds for any $t\in [0,T]$ while the sum of the RHS terms is non-decreasing in $t$, we then see
\begin{eqnarray*}
\sup_{0\le u\le t} \E|x_\D(u)|^{\bp}
\le C + C \int_0^t \Big( \sup_{0\le u\le s} \E|x_\D(u)|^{\bp} \Big)ds.
\end{eqnarray*}
The well-known Gronwall inequality yields that
$$
\sup_{0\le u\le T} \E|x_\D(u)|^{\bp} \le C.
$$
As this holds for any $\D \in (0,\D^*]$ while $C$ is independent of $\D$, we see the
required assertion (\ref{4.4}). $\Box$

The following two lemmas are the analogues of Lemmas \ref{L3.2} and \ref{L3.3}.

\begin{lemma}\label{L4.5A}
Let Assumptions \ref{A2.1} and \ref{A4.1} hold. For any real number $R> \|\xi\|$,
define the stopping time
$\T_R=\inf\{t\ge 0: |x(t)|\ge R\}$.
Then
\begin{equation}\label{4.6a}
\PP(\T_R\le T) \le \frac{C}{ R^{{}^{\bp}}}.
\end{equation}
\end{lemma}

\begin{lemma}\label{L4.5B}
Let Assumptions \ref{A2.1} and \ref{A4.1} hold. For any real number $R> \|\xi\|$
and $\D\in (0,\D^*]$, define the stopping time
$\r_{\D,R} =\inf\{t\ge 0: |x_\D(t)|\ge R\}$.
Then
\begin{equation}\label{5.7b}
\PP(\r_{\D,R}\le T) \le \frac{C}{R^{{}^{\bp}}}.
\end{equation}
\end{lemma}

Their proofs are similar to those of  Lemmas \ref{L3.2} and \ref{L3.3}, respectively,
so are omitted.  We can now state our main result in this section.

\begin{theorem}\label{T4.5}
Let Assumptions \ref{A2.1}, \ref{A5.2} and \ref{A4.1} hold.  Then, for any $q\in [2,\bp)$,
\begin{equation}\label{4.6}
\lim_{\D\to 0} \E|x_\D(T) - x(T)|^q = 0
\quad\hbox{and}\quad
\lim_{\D\to 0} \E|\bar x_\D(T) - x(T)|^q = 0.
\end{equation}
\end{theorem}

\noindent
{\it Proof}.   We use the same notation as in the proof of Theorem \ref{T3.4}.
Fix any $q \in [2,\bp)$.
Using the Young inequality, we can show that for any $\de > 0$,
\begin{eqnarray} \label{4.6b}
\E|e_\D(T)|^q
&\le & \E\Big( |e_\D(T)|^q I_{\{\o_{\D,R}>T\}} \Big)
+ \frac{q\de}{\bp} \E|e_\D(T)|^{\bp} + \frac{\bp-q}{\bp\de^{q/(\bp-q)}} \PP(\o_{\D,R}\le T).
\end{eqnarray}
By Lemmas \ref{L4.3} and \ref{L4.4}, we have
\begin{eqnarray} \label{4.6c}
\E|e_\D(T)|^{\bp} \le C,
\end{eqnarray}
while by Lemmas \ref{L4.5A} and \ref{L4.5B},
\begin{eqnarray} \label{4.6d}
\PP(\o_{\D,R}\le T) \le \PP(\T_R\le T)+\PP(\r_{\D,R}\le T) \le \frac{C}{R^{{}^{\bp}}}.
\end{eqnarray}
Using these and (\ref{3.20}) (please recall Remark \ref{R3.5}),
we obtain
\begin{eqnarray} \label{4.7}
\E|e_\D(T)|^q \le H\D^{q(0.5\we\g)}
+ \frac{C q\de }{\bp} + \frac{C(\bp-q)}{\bp R^{{}^{\bp}} \de^{q/(\bp-q)}}.
\end{eqnarray}
Now, for any $\e>0$, we first choose $\de$ sufficiently small for
$Cq\de/\bp \le \e/3$ and then choose $R$ sufficiently large for
$$
\frac{C(\bp-q)}{\bp R^{{}^{\bp}} \de^{q/(\bp-q)}} \le \frac{\e}{3},
$$
and further then choose $\D$ sufficiently small for $H\D^{q(0.5\we\g)} \le \e/3$ to get that
\begin{equation} \label{4.9}
\E|e_\D(T)|^q \le \e.
\end{equation}
In other words, we have shown that
$$
\lim_{\D\to 0} \E|e_\D(T)|^q = 0.
$$
This, along with Lemma \ref{L2.5}, implies another assertion
$$
\lim_{\D\to 0} \E|x(T)-\bar x_\D(T)|^q = 0.
$$
The proof is therefore complete. $\Box$

Let us now discuss an example to illustrate this theorem before we
study the convergence rates.

\begin{expl} \label{E5.5}
{\rm Consider the scalar SDDE
\begin{equation} \label{E5.5a}
dx(t) = f(x(t),x(t-\T))dt + g(x(t),x(t-\T))dB(t), \quad t\ge 0,
\end{equation}
with the initial data $\{x(\o):-\T\le\o\le 0\}=\xi\in C([-\T,0];\RR)$ which
satisfy Assumption \ref{A5.2}, where
$$
f(x,y) = a_1 +a_2|y|^{4/3} - a_3x^3  \quad\hbox{and}\quad
g(x,y) = a_4|x|^{3/2}+ a_5y, \quad x,y\in\RR,
$$
and $a_1,\cdots,a_5$ are all real numbers with $a_3>0$.  Clearly, the coefficients
$f$ and $g$
are locally Lipschitz continuous, namely, they satisfy Assumption \ref{A2.1}.
Moreover, for any $\bp >2$, we have
\begin{eqnarray*}
xf(x,y)+\frac{\bp-1}{2} |g(x,y)|^2
&\le& |a_1||x|+|a_2||x||y|^{4/3}-a_3|x|^4+ (\bp-1)(|a_4||x|^3+|a_5||y|^2).
\end{eqnarray*}
But, by the Young inequality,
$$
|x||y|^{4/3} = (|x|^3)^{1/3} (|y|^2)^{2/3} \le  |x|^3+|y|^2.
$$
We therefore have
\begin{eqnarray*}
& & xf(x,y)+\frac{\bp-1}{2} |g(x,y)|^2  \\
& \le & |a_1||x|+(|a_2|+|a_4|(\bp-1))|x|^3  -a_3|x|^4+ (|a_2|+a_5(\bp-1))|y|^2 \\
&\le &  K_1 (1+ |y|^2),
\end{eqnarray*}
where $K_1=(|a_2|+|a_5|(\bp-1))\ve K$ and
$$
K = \sup_{u\ge 0} \big[|a_1|u+(|a_2|+|a_4|(\bp-1))u^3  -a_3u^4\big] <\8.
$$
That is, Assumption \ref{A4.1} is satisfied for any $\bp >2$.
To apply Theorem \ref{T4.5}, we still need to design
functions $\mu$ and $h$ satisfying (\ref{mudef}) and (\ref{hdef}).
Note that
$$
\sup_{|x|\le u} (|f(x)|\ve|g(x)|) \le \hat a u^3, \quad\forall u\ge 1,
$$
where $\hat a= (|a_1|+|a_2|+a_3)\ve (|a_4|+|a_5|)$.
We can hence have $\mu(u)=\hat a u^3$ and its inverse function
$\mu^{-1}(u)=(u/\hat a)^{1/3}$ for $u\ge 0$.  For $\e \in (0,1/4]$, we define
$h(\D)=\D^{-\e}$ for $\D>0$.  Letting $\D^*\in (0,1]$ be sufficiently small,
we can make (\ref{hdef}) hold.
By Theorem \ref{T4.5}, we can then conclude  that the truncated EM solutions
will converge to the true solution $x(t)$  in the sense that
$$
\lim_{\D\to 0} \E|x_\D(T) - x(T)|^q = 0
\quad\hbox{and}\quad
\lim_{\D\to 0} \E|\bar x_\D(T) - x(T)|^q = 0
$$
for every $q\ge 2$.
}
\end{expl}

\section{Convergence Rates}\label{sec:5}

In the previous sections, we showed the convergence in $L^q$ of the
truncated EM solutions to the true solution. However, the convergence
was in the asymptotic form without the convergence rate.
In this section we will discuss the rate.  To avoid the notation becoming too complicated,
we will only discuss the convergence rate in $L^2$ but the technique developed here can
certainly be applied to study the rate in $L^q$.
Recall that we use two functions  $\mu(\cdot)$ and $h(\cdot)$ to define the
truncated EM method.  The choices of these functions are independent as long as they satisfy (\ref{mudef}) and (\ref{hdef}), respectively.
It is interesting to see that they will satisfy a related condition in order for us to obtain the convergence rate.

We need an additional condition.  To state it, we need a new notation.
Let ${\cal U}$ denote the family of continuous functions $U:\RR^n\K\RR^n \to \RR_+$
such that for each $b>0$, there is a positive constant $\k_b$ for which
$$
U(x,\bx)\le \k_b|x-\bx|^2,\ \ \ \forall x,\bx\in\RR^n \hbox{ with } |x|\ve |\bx|\le b.
$$

\begin{assp}\label{A5.1}
Assume that there is a positive constant $H_1$ and a  function $U \in {\cal U}$ such that
\begin{eqnarray}\label{5.1}
& & (x-\bx)^T(f(x,y)-f(\bx,\bar y)) + \frac{1}{2} |g(x,y)-g(\bx,\bar y)|^2
\nonumber \\
& & \ \ \le
H_1(|x-\bx|^2+|y-\bar y|^2) -U(x,\bx)+U(y,\bar y)
\end{eqnarray}
for all $x,y,\bx, \bar y\in\RR^n$.
\end{assp}

Let us first present a key lemma.

\begin{lemma}\label{L5.3}
Let Assumptions \ref{A2.1}, \ref{A5.2} and \ref{A5.1} hold.
Let $R>\|\xi\|$ be a real number and let $\D\in (0,\D^*)$ be sufficiently small
such that $\mu^{-1}(h(\D)) \ge R$. Let $\o_{\D,R}$ and $e_\D(t)$
be the same as defined in Section 3. Then
\begin{equation}\label{5.3}
\E|e_\D(T\we \o_{\D,R})|^2 \le C (\D^{2\g} \ve [\D^{1/2} (h(\D))^2]),
\end{equation}
where, as before, $C$ is the generic constant independent of $R$ and $\D$.
\end{lemma}

\noindent
{\it Proof}. We write $\o_{\D,R}=\o$ for simplicity.
The It\^o formula shows that
\begin{eqnarray}\label{5.4a}
\E|e_\D(t\we\o)|^2 & = &
\E\int_0^{t\we\o} \Big( 2 e_\D^T(s) [f(x(s),x(s-\T))-f_\D(\bx_\D(s),\bx_\D(s-\T))]
\nonumber \\
& & \quad\quad +
|g(x(s),x(s-\T))-g_\D(\bx_\D(s),\bx_\D(s-\T))|^2 \Big) ds
\end{eqnarray}
for $0\le t \le T$. We observe that for $0\le s\le t\we\o$,
$$
|\bar x_\D(s)| \ve |\bar x_\D(s-\T)|\ve |x(s)| \ve |x(s-\T)|\le R.
$$
But we have the condition that $\mu^{-1}(h(\D)) \ge R$, so
$$
|\bar x_\D(s)| \ve |\bar x_\D(s-\T)|\ve |x(s)| \ve |x(s-\T)|\le \mu^{-1}(h(\D)).
$$
Recalling the definition of the truncated
functions $f_\D$ and $g_\D$ as well as (\ref{mudef}), we hence have that
$$
f_\D(\bx_\D(s),\bx_\D(s-\T))=f(\bx_\D(s),\bx_\D(s-\T)),
\ \
g_\D(\bx_\D(s),\bx_\D(s-\T))=g(\bx_\D(s),\bx_\D(s-\T))
$$
and
\begin{equation} \label{5.4}
|f(x(s),x(s-\T))|\ve |f(\bx_\D(s),\bx_\D(s-\T))|\le h(\D)
\end{equation}
for $0\le s\le t\we\o$.  It therefore follows from (\ref{5.4a}) that
\begin{eqnarray}\label{5.4c}
& & \E|e_\D(t\we\o)|^2 \nonumber  \\
& =& \E\int_0^{t\we\o} \Big( 2 e_\D^T(s) [f(x(s),x(s-\T))-f(\bx_\D(s),\bx_\D(s-\T))]
\nonumber  \\
& & \quad\quad +
|g(x(s),x(s-\T))-g(\bx_\D(s),\bx_\D(s-\T))|^2 \Big) ds   \\
& =& \E\int_0^{t\we\o} \Big( 2 (x(s)-\bx_\D(s))^T [f(x(s),x(s-\T))-f(\bx_\D(s),\bx_\D(s-\T))]
 \nonumber  \\
& & \quad\quad +
|g(x(s),x(s-\T))-g(\bx_\D(s),\bx_\D(s-\T))|^2 \Big) ds \nonumber \\
& + &  \E\int_0^{t\we\o}  2 (\bx_\D(s)-x_\D(s))^T [f(x(s),x(s-\T))-f(\bx_\D(s),\bx_\D(s-\T))]ds. \nonumber
\end{eqnarray}
By Assumption \ref{A5.1} and (\ref{5.4}),   we then derive that
\begin{eqnarray} \label{5.5}
\E|e_\D(t\we\o)|^2
&\le&  2H_1 \E\int_0^{t\we\o} \Big(|x(s)-\bar x_\D(s)|^2+ |x(s-\T)-\bar x_\D(s-\T)|^2 \Big) ds
\nonumber \\
&+& \E\int_0^{t\we\o} \Big( -U(x(s),\bar x_\D(s)) + U(x(s-\T),\bar x_\D(s-\T)) \Big)ds
\nonumber  \\
& +& 4h(\D) \E\int_0^{t\we\o}  |\bar x_\D(s)-x_\D(s)| ds.
\end{eqnarray}
But, by Assumption \ref{A5.2} and Lemma \ref{L2.5}, we derive that
\begin{eqnarray} \label{5.6}
 & & \E\int_0^{t\we\o} \Big(|x(s)-\bar x_\D(s)|^2+ |x(s-\T)-\bar x_\D(s-\T)|^2 \Big) ds
\nonumber \\
&\le & 2\E\int_0^{t\we\o} \Big(|e_\D(s)|^2+|e_\D(s-\T)|^2
+|x_\D(s)-\bar x_\D(s)|^2+ |x_\D(s-\T)-\bar x_\D(s-\T)|^2 \Big) ds
\nonumber  \\
&\le & 4\E\int_0^{t}  |e_\D(s\we\o)|^2 ds
+ 4 \int_0^{T} \E|x_\D(s)-\bar x_\D(s)|^2ds +\int_{-\T}^0 |\xi(s)- \xi(\lfloor s/\D\rfloor \D)|^2 ds
\nonumber  \\
&\le & 4 \int_0^{t}  \E|e_\D(s\we\o)|^2 ds
+ C\D (h(\D))^2 +  \T K_3^2 \D^{2\g}.
\end{eqnarray}
Moreover, by the property of the ${\cal U}$-class function $U$ and Assumption \ref{A5.2},
we have
\begin{eqnarray} \label{5.7}
& & \E\int_0^{t\we\o} \Big( -U(x(s),\bar x_\D(s)) + U(x(s-\T),\bar x_\D(s-\T)) \Big)ds
\nonumber \\
&\le&  \int_{-\T}^0 U(\xi(s),\xi(\lfloor s/\D\rfloor \D)) ds
 \le \int_{-\T}^0 \k_b |\xi(s)-\xi(\lfloor s/\D\rfloor \D)|^2 ds \nonumber\\
& \le &  \T \k_b K_3^2  \D^{2\g},
\end{eqnarray}
where $b=\|\xi\|$. Furthermore, by Lemma \ref{L2.5},
\begin{equation}\label{5.8}
\E\int_0^{t\we\o}  |\bar x_\D(s)-x_\D(s)| ds \le
 \int_0^{T}  \E|\bar x_\D(s)-x_\D(s)| ds \le C\D^{1/2}h(\D).
\end{equation}
Substituting (\ref{5.6})-(\ref{5.8}) into (\ref{5.5}), we get
 \begin{eqnarray*}
\E|e_\D(t\we\o)|^2 \le
8H_1 \int_0^{t} \E|e_\D(s\we\o)|^2 ds +C (\D^{2\g} \ve [\D^{1/2} (h(\D))^2]).
\end{eqnarray*}
By the Gronwall inequality, we obtain the required assertion (\ref{5.3}).
$\Box$

Let us now state our first result on the convergence rate, where we reveal a strong relation between
functions $\mu(\cdot)$ and $h(\cdot)$, which are used to define the
truncated EM method.

\begin{theorem}\label{T5.4}
Let Assumptions \ref{A2.1},  \ref{A5.1}, \ref{A4.1} and \ref{A5.2} hold.
Assume that
\begin{equation}\label{5.8b}
h(\D) \ge \mu\big( (\D^{2\g} \ve [\D^{1/2} (h(\D))^2])^{-1/(\bp-2)} \big)
\end{equation}
for all sufficiently small $\D\in (0,\D^*)$.  Then, for every such small $\D$,
\begin{equation}\label{5.9}
\E|x(T)-x_\D(T)|^2 \le C (\D^{2\g} \ve [\D^{1/2} (h(\D))^2])
\end{equation}
and
\begin{equation}\label{5.10}
\E|x(T)-\bar x_\D(T)|^2 \le C (\D^{2\g} \ve [\D^{1/2} (h(\D))^2]).
\end{equation}
\end{theorem}

\noindent
{\it Proof}. We use the same notation as in the proof of Theorem \ref{T4.5}.
It follows from (\ref{4.6b})-(\ref{4.6d}) with $q=2$ that the inequality
\begin{eqnarray} \label{5.11}
\E|e_\D(T)|^2
&\le & \E\Big( |e_\D(T\we\o_{\D,R})|^2 \Big)
+\frac{2C\de}{\bp}
+   \frac{C(\bp-2)}{\bp R^{{}^{\bp}} \de^{2/(\bp-2)}}
\end{eqnarray}
holds for any $\D\in (0,\D^*)$, $R> \|\xi\|$ and $\de>0$.   In
particular,  choosing
$$
\de = \D^{2\g} \ve [\D^{1/2} (h(\D))^2]\quad\hbox{and}\quad
R = (\D^{2\g} \ve [\D^{1/2} (h(\D))^2])^{-1/(\bp-2)},
$$
we get
\begin{equation} \label{5.13}
\E|e_\D(T)|^2 \le \E|e_\D(T\we \o_{\D,R})|^2 +  C (\D^{2\g} \ve [\D^{1/2} (h(\D))^2]).
\end{equation}
But, by condition (\ref{5.8b}), we have
$$
\mu^{-1}(h(\D)) \ge(\D^{2\g} \ve [\D^{1/2} (h(\D))^2])^{-1/(\bp-2)} = R.
$$
We can hence apply Lemma \ref{L5.3} to obtain
\begin{equation} \label{5.14}
\E|e_\D(T\we \o_{\D,R})|^2 \le C (\D^{2\g} \ve [\D^{1/2} (h(\D))^2]).
\end{equation}
Substituting this into (\ref{5.13}) yields  the first assertion (\ref{5.9}) . The second
assertion (\ref{5.10})
follows from (\ref{5.9}) and Lemma \ref{L2.5}.  $\Box$

Let us discuss an example to illustrate Theorem \ref{T5.4} and to motivate
our further results on the convergence rates.

\begin{expl} \label{E5.5B}
{\rm
 Consider the same SDE in Example \ref{E5.5}.
 We need to verify Assumption \ref{A5.1}.
For $x,y,\bx,\bar y\in \RR$, it is easy to show that
\begin{equation}\label{E5.5b}
  (x-\bx)(f(x,y)-f(\bx, \bar y))\le a_2^2|x-\bx|^2 +(|y|^{4/3}-|\bar y|^{4/3})^2
-0.5a_3|x-\bx|^2 (x^2+\bx^2).
\end{equation}
But, by the mean value theorem,
$$
 (|y|^{4/3}-|\bar y|^{4/3})^2 \le \frac{16}{9}  |y-\bar y|^2   (|y|^{1/3}+|\bar y|^{1/3})^2
 \le 4 |y-\bar y|^2   (|y|^{2/3}+|\bar y|^{2/3}).
$$
Let $a_6:=\sup_{u\ge 0} (8u^{2/3}-0.5a_3 u^2)$.  Then $0\le a_6<\8$ and
$$
 (|y|^{4/3}-|\bar y|^{4/3})^2 \le  a_6|y-\bar y|^2 + 0.25a_3|y-\bar y|^2(y^2+\bar y^2).
 $$
 Substituting this into (\ref{E5.5b}) yields
 \begin{eqnarray} \label{E5.5c}
  & & (x-\bx)(f(x,y)-f(\bx, \bar y))\nonumber \\
&  \le &  (a_6\ve a_2^2)(|x-\bx|^2+|y-\bar y|^2) \nonumber \\
& -  & 0.5a_3|x-\bx|^2 (x^2+\bx^2) + 0.25a_3|y-\bar y|^2(y^2+\bar y^2).
\end{eqnarray}
Similarly, we can show that
\begin{equation}\label{E5.5d}
 0.5|g(x,y)-g(\bx, \bar y) |^2
\le (a_7\ve a_5^2)(|x-\bx|^2+|y-\bar y|^2)
+ 0.25a_3|x-\bx|^2 (x^2+\bx^2),
\end{equation}
where $a_7:=\sup_{u\ge 0}(9a_4^2u-0.5a_3u^2) \in (0,\8)$.
It then follows from (\ref{E5.5c}) and (\ref{E5.5d}) that
\begin{eqnarray}\label{E5.5e}
  & & (x-\bx)(f(x,y)-f(\bx, \bar y))  + 0.5|g(x,y)-g(\bx, \bar y) |^2 \nonumber \\
 & & \   \le   H_1(|x-\bx|^2+|y-\bar y|^2) - U(x,\bx)+U(y,\bar y),
\end{eqnarray}
where $H_1= (a_6\ve a_2^2)+ (a_7\ve a_5^2)$ and
$U(x,\bx)= 0.25a_3|x-\bx|^2 (x^2+\bx^2)$.  It is obvious that $U\in {\cal U}$.
In other words, we have shown that Assumption \ref{A5.1} is satisfied too.
To apply Theorem \ref{T5.4}, we use the same functions $\mu(\cdot)$ and $h(\cdot)$
as defined in Example \ref{E5.5}.
We observe that inequality (\ref{5.8b}) becomes
\begin{equation}\label{5.20}
\D^{-\e} \ge \hat a \D^{ -3[(2\g)\we (1/2-2\e)]/(\bp-2)}.
\end{equation}
But, for any $\e \in (0,1/4]$, we can choose $\bp$ sufficiently large such that
$\e > 3[(2\g)\we (1/2-2\e)]/(\bp-2)$ and hence (\ref{5.20}) holds for all
sufficiently small $\D$.  We can therefore conclude by Theorem \ref{T5.4} that
the truncated EM solutions of the SDE (\ref{E5.5a}) satisfy
\begin{equation}\label{5.21}
\E|x(T)-x_\D(T)|^2 = O( \D^{(2\g)\we (1/2-2\e)})
\  \ \hbox{and}\ \
\E|x(T)-\bar x_\D(T)|^2= O( \D^{(2\g)\we (1/2-2\e)}).
\end{equation}
}
\end{expl}

It is known that for every $\a\in (0,0.5)$, the Brownian motion is $\a$-H\"older continuous
(see, e.g., \cite{KR88}).  If we regard the initial data $\xi(u)$, $u\in [-\T,0]$ as an observation
of the state during the time interval $[-\T,0]$, it is reasonable to assume that $\g\in (0,0.5)$.
If $\g$ is close to 0.5, then (\ref{5.21}) shows the order of  convergence is close to 0.25.  Can we improve the order?  The answer is yes though we need stronger conditions.

\begin{assp}\label{A5.6}
Assume that there are positive constants $\a$ and $H_2$ and a  function $U \in {\cal U}$ such that
\begin{eqnarray}\label{5.23}
& & (x-\bx)^T(f(x,y)-f(\bx,\bar y)) + \frac{1+\a}{2} |g(x,y)-g(\bx,\bar y)|^2
\nonumber \\
& & \ \ \le
H_2(|x-\bx|^2+|y-\bar y|^2) -U(x,\bx)+U(y,\bar y)
\end{eqnarray}
for all $x,y,\bx, \bar y\in\RR^n$.
\end{assp}

\begin{assp}\label{A5.7}
Assume that there is a pair of positive constants $r$ and $H_3$   such that
\begin{eqnarray}\label{5.24}
& & |f(x,y)-f(\bx,\by)|^2\ve |g(x,y)-g(\bx,\by)|^2
\nonumber \\
&\le &
H_3(|x-\bx|^2+|y-\bar y|^2)(1+|x|^r+|\bx|^r+|y|^r+|\by|^r)
\end{eqnarray}
for all $x,y,\bx, \bar y\in\RR^n$.
\end{assp}

\begin{lemma}\label{L5.8}
Let Assumptions \ref{A2.1}, \ref{A5.2}, \ref{A4.1}, \ref{A5.6} and \ref{A5.7} hold
and $\bp > r$.
Let $R>\|\xi\|$ be a real number and let $\D\in (0,\D^*)$ be sufficiently small
such that $\mu^{-1}(h(\D)) \ge R$.  Let $\o_{\D,R}$ and $e_\D(t)$
be the same as defined in Section 3. Then
\begin{equation}\label{5.25}
\E|e_\D(T\we \o_{\D,R})|^2 \le C (\D^{2\g} \ve [\D (h(\D))^2]).
\end{equation}
\end{lemma}

\noindent
{\it Proof}.  We use the same notation as in the proof of Lemma \ref{L5.3}.
It follows from (\ref{5.4c}) that
\begin{eqnarray}\label{5.26}
\E|e_\D(t\we\o)|^2
& \le & \E\int_0^{t\we\o} \Big( 2 e_\D^T(s) [f(x(s),x(s-\T))-f(x_\D(s),x_\D(s-\T))] \nonumber \\
& +&
(1+\a) |g(x(s),x(s-\T))-g(x_\D(s),x_\D(s-\T))|^2 \nonumber \\
& +& 2 e_\D^T(s) [f(x_\D(s),x_\D(s-\T))-f(\bx_\D(s),\bx_\D(s-\T))] \nonumber  \\
& +&  (1+\a^{-1})|g(x_\D(s),x_\D(s-\T))-g(\bx_\D(s),\bx_\D(s-\T))|^2 \Big) ds.
 \end{eqnarray}
By Assumptions \ref{A5.2}, \ref{A5.6} and \ref{A5.7}, we can then show
\begin{equation}\label{5.27}
\E|e_\D(t\we\o)|^2
 \le  (4H_2+1)  \int_0^t  \E|e_\D(s\we\o)|^2 ds + 2\T\k_bK_3^2\D^{2\g} + J,
\end{equation}
where (\ref{5.7}) has been used and
\begin{eqnarray*}
J & :=  &  \E\int_0^{t\we\o}
H_3(2+\a^{-1}) (|x_\D(s)-\bx_\D(s)|^2+|x_\D(s-\T)-\bx_\D(s-\T)|^2)    \\
& & \quad\quad \K (1+|x_\D(s)|^r+|\bx_\D(s)|^r+|x_\D(s-\T)|^r+|\bx_\D(s-\T)|^r) ds.
 \end{eqnarray*}
But, by the H\"older inequality, Lemmas \ref{L2.5} and \ref{L4.3} and Assumption \ref{A5.2}, we can derive that
\begin{eqnarray*}
J & \le  & C   \int_0^T
\Big(\E|x_\D(s)-\bx_\D(s)|^{2\bp/(\bp-r)}
+\E|x_\D(s-\T)-\bx_\D(s-\T)|^{2\bp/(\bp-r)}  \Big)^{(\bp-r)/\bp}    \\
& & \quad\quad \K \Big(  1+\E|x_\D(s)|^{\bp}+\E|\bx_\D(s)|^{\bp}+\E|x_\D(s-\T)|^{\bp}
+ \E|\bx_\D(s-\T)|^{\bp}   \Big)^{r/\bp} ds \\
&\le & C(\D^{2\g}\ve [\D(h(\D))^2]).
 \end{eqnarray*}
Substituting this into (\ref{5.27}) gives
\begin{equation*}
\E|e_\D(t\we\o)|^2
 \le  (4H_2+1)  \int_0^t  \E|e_\D(s\we\o)|^2 ds + C(\D^{2\g}\ve [\D(h(\D))^2]),
\end{equation*}
which implies the required assertion (\ref{5.25}).  $\Box$

The following theorem gives a better convergence rate than Theorem \ref{T5.4}.

\begin{theorem}\label{T5.9}
Let Assumptions \ref{A2.1}, \ref{A5.2}, \ref{A4.1},  \ref{A5.6}  and \ref{A5.7} hold
and $\bp >r$.
Assume that
\begin{equation}\label{5.28}
h(\D) \ge \mu\big( (\D^{2\g} \ve [\D(h(\D))^2])^{-1/(\bp-2)} \big)
\end{equation}
for all sufficiently small $\D\in (0,\D^*)$.  Then, for every such small $\D$,
\begin{equation}\label{5.29}
\E|x(T)-x_\D(T)|^2 \le C (\D^{2\g} \ve [\D(h(\D))^2])
\end{equation}
and
\begin{equation}\label{5.30}
\E|x(T)-\bar x_\D(T)|^2 \le C (\D^{2\g} \ve [\D(h(\D))^2]).
\end{equation}
\end{theorem}

\noindent
{\it Proof}. We use the same notation as in the proof of Theorem \ref{T5.4}.
Choosing
$$
\de = \D^{2\g} \ve [\D(h(\D))^2]\quad\hbox{and}\quad
R = (\D^{2\g} \ve [\D(h(\D))^2])^{-1/(\bp-2)},
$$
we get from (\ref{5.11}) that
\begin{equation} \label{5.31}
\E|e_\D(T)|^2 \le \E|e_\D(T\we \o_{\D,R})|^2 +  C (\D^{2\g} \ve [\D (h(\D))^2]).
\end{equation}
But, by condition (\ref{5.28}), we have
$$
\mu^{-1}(h(\D)) \ge(\D^{2\g} \ve [\D (h(\D))^2])^{-1/(\bp-2)} = R.
$$
We can hence apply Lemma \ref{L5.8} to obtain
\begin{equation} \label{5.32}
\E|e_\D(T\we \o_{\D,R})|^2 \le C (\D^{2\g} \ve [\D (h(\D))^2]).
\end{equation}
Substituting this into (\ref{5.31}) yields  the first assertion (\ref{5.29}) . The second
assertion (\ref{5.30})
follows from (\ref{5.29}) and Lemma \ref{L2.5}.  $\Box$

\begin{expl}\label{E5.5C}
{\rm  Let us return to Example \ref{E5.5}
once again.
Instead of (\ref{E5.5d}), we can have the following alternative estimate
\begin{equation}\label{E5.5g}
 |g(x,y)-g(\bx, \bar y) |^2
\le 2(a_8\ve a_5^2)(|x-\bx|^2+|y-\bar y|^2)
+ 0.25a_3|x-\bx|^2 (x^2+\bx^2),
\end{equation}
where $a_8:=\sup_{u\ge 0}(9a_4^2u-0.25a_3u^2) \in (0,\8)$.    It then follows from (\ref{E5.5c}) and (\ref{E5.5g}) that
\begin{eqnarray}\label{E5.5h}
  & & (x-\bx)(f(x,y)-f(\bx, \bar y))  +  |g(x,y)-g(\bx, \bar y) |^2\nonumber \\
 & & \   \le   H_2(|x-\bx|^2+|y-\bar y|^2) - U(x,\bx)+U(y,\bar y),
\end{eqnarray}
where $H_2= (a_6\ve a_2^2)+ 2(a_8\ve a_5^2)$ and
$U(x,\bx)= 0.25a_3|x-\bx|^2 (x^2+\bx^2)$.
In other words, we have shown that Assumption \ref{A5.6} is satisfied with $\a=1$.  It is
also straightforward to show that
\begin{equation}\label{E5.5i}
|f(x,y)-f(\bx,\by)|^2 \le 8a_2^2|y-\by|^2(1+|y|^4+|\by|^4)+16a_3^2|x-\bx|^2(|x|^4+|\bx|^4).
\end{equation}
We hence see from (\ref{E5.5g}) and (\ref{E5.5i})  that Assumption \ref{A5.7}
is also satisfied with $r=4$.   In other words, we have shown that
Assumptions \ref{A2.1}, \ref{A4.1}, \ref{A5.2}, \ref{A5.6}  and \ref{A5.7} hold
for every $\bp >r=4$.
Let $\mu(\cdot)$ and
$h(\cdot)$ be the same as before.  We can then conclude by Theorem \ref{T5.9} that
the truncated EM solutions of the SDE (\ref{E5.5a}) satisfy
\begin{equation}\label{E5.5j}
\E|x(T)-x_\D(T)|^2 = O( \D^{(2\g)\we (1-2\e)})
\  \ \hbox{and}\ \
\E|x(T)-\bar x_\D(T)|^2= O( \D^{(2\g)\we (1-2\e)}).
\end{equation}
In particular, if $\g$ is close to 0.5 (or bigger than half), this shows that the order of convergence is close to 0.5.
}
\end{expl}

\section{Conclusion}
In this paper we have used the
new explicit method, called the truncated EM method, to study the
strong convergence of the numerical solutions for nonlinear SDDEs.
For a given stepsize $\D$, we define the discrete-time truncated EM numerical solution and then form
two versions of the continuous-time truncated EM solutions, namely the continuous-time
step-process truncated EM solution $\bar x_\D(t)$ and
the continuous-time continuous-process truncated EM solution $x_\D(t)$.
Under the local Lipschitz condition plus the generalized
Khasminskii-type condition, we have successfully shown the strong convergence
of both continuous-time truncated EM solutions to the true solution in the sense that
$$
\lim_{\D\to 0} \E|x_\D(T) - x(T)|^q = 0
\quad\hbox{and}\quad
\lim_{\D\to 0} \E|\bar x_\D(T) - x(T)|^q = 0
$$
for any $T >0$ and  $q\in [1,2)$. Under a slightly stronger Khasminskii-type condition,
we have showed the above convergence for some $q\ge 2$.
We have also discussed the convergence rates in $L^2$ under some additional conditions.
We have used several examples to illustrate our theory throughout the paper.

\section*{Acknowledgements}

The authors would
like to thank the Leverhulme Trust (RF-2015-385),
the Royal Society of London (IE131408),
%the EPSRC (EP/E009409/1),
%the Royal Society of Edinburgh (RKES115071),
%the London Mathematical Society (11219),
%the Edinburgh Mathematical Society (RKES130172),
the Natural Science Foundation of China (11471216), %Yue
the Natural Science Foundation of Shanghai (14ZR1431300), %Guo
the Ministry of Education (MOE) of China (MS2014DHDX020)
for their financial support.

\end{document}